\documentclass[a4paper,12pt]{article}
\usepackage[english]{babel}
\usepackage{amsmath}
\usepackage{amssymb}
\usepackage{tabls}
\usepackage[pdftex]{graphicx}
\usepackage{calc}
\usepackage{xcolor}
\usepackage[T1]{fontenc}
\usepackage[latin1]{inputenc}
\title{The Center of Crystalline Graded Rings}
\author{Tim Neijens\\University of Antwerp\\ \texttt{tim.neijens@gmail.com} \and Freddy Van Oystaeyen \\ University of Antwerp \\ \texttt{fred.vanoystaeyen@ua.ac.be}}

\newcommand{\Z}{\mathbb{Z}}
\newcommand{\Q}{\mathbb{Q}}

\newcommand{\blok}{\hfill \Box}
\newcommand{\CGR}{\mathop \diamondsuit \limits_{\sigma ,\alpha}}
\newtheorem{defi}{Definition}[section]
\newtheorem{opm}[defi]{Remark}
\newtheorem{prop}[defi]{Proposition}
\newtheorem{gev}[defi]{Corollary}
\newtheorem{st}[defi]{Theorem}
\newtheorem{lem}[defi]{Lemma}

\begin{document}
\maketitle

\begin{abstract}
In the first section of the paper, we will give some basic definitions and properties about Crystalline Graded Rings.  In the following section we will provide a general description of the center.  Afterwards, the case where the grading group is Abelian finite will be handled.  The center will have some properties of a crystalline graded ring, but not all.  We will call this Arithmetically Crystalline Graded.  The center is crystalline graded if the part of degree zero is a principal ideal domain.  The last section deals with the case where the grading group is non-Abelian finite.  Since this situation is much more complicated than the Abelian case, we primarily focus on the conditions to have a trivial center.  The fact that the center is Arithmetically Crystalline Graded also holds in this case.
\end{abstract}

\section*{Introduction}
Crystalline graded rings have been introduced in \cite{NVO6} as a generalization of crossed products on one hand and of generalized Weyl algebras on the other, cf. \cite{B1}, \cite{BVO3}.  In this paper we focus on the situation where the subring of degree $e$, say $A_e$, of the crystalline graded ring $A$, is a Dedekind domain, $A_e = D$ with field of fractions $L$.  Important examples in case the grading group is torsionfree Abelian, usually $\Z$ or $\Z^n$, include $D$-orders in well-known algebras over $L$ like the first Weyl algebra, the quantum plane, quantum $sl_2$, \ldots, many algebras of quantum type appearing in : \cite{B}, \cite{BVO3}, \cite{NVO6}, \ldots
However when attention is restricted to finite grading groups $G$, then crystalline graded rings over a commutative ring are necessarily PI-rings, then typical examples of quantum-type (quantized or $q$-deformed algebras) would correspond to the case where the deformation parameter $q$ is a root of unity.  Hence the class of crystalline graded rings over $D$ consists of $D$-orders in classical crossed products over the field of fractions which are finite dimensional $L$-algebras.  It is a natural problem to aim at an algebraic classification for crystalline graded rings over a Dedekind domain $D$, specializing to particular cases like discrete valuation rings, $k[T]$ or rings of integers in number fields for more concrete applications.
On the other hand, most results extend to the situation where $A_e$ is the coordinate ring of a normal variety of dimension $d$, i.e. a Noetherian integrally closed domain (having global dimension $d$), but we do not consider this extension here.  Applying the theory of maximal $D$-orders seems to be the obvious way to start analyzing the algebra structure of the noncommutative algebras under consideration but at an even more elementary level there is first the problem of determining the center of $A$.
The first part of this paper deals with the determination of the center of a crystalline graded ring as before.  Knowing that skew group rings and twisted group rings are special cases of our constructions, it is clear that rings of invariants for certain specific group actions and ray classes, as in the theory of projective representations of finite groups, will play an important part in this.  Secondly, the particular case where the center $Z(A)$ is minimal (i.e. $Z(A)=D$ or $Z(A)= D^G$) should be the most easy to describe.  In fact, the results concerning Clifford representations for Abelian groups, cf. \cite{CVO}, suggest that such crystalline graded rings will be (maximal) $D$-orders in generalized Clifford algebras over the quotient field $L$ of $D$.

\section{Preliminaries}

{\defi \label{def1}\textbf{Pre-Crystalline Graded Ring}\\
Let $A$ be an associative ring with unit $1_A$.  Let $G$ be an arbitrary group.  Consider an injection $u: G \rightarrow A$ with $u_e = 1_A$, where $e$ is the neutral element of $G$ and $u_g \neq 0$,  $\forall g \in G$.  Let $R \subset A$ be an associative ring with $1_{R}=1_A$.  We consider the following properties:
\begin{description} 
	\item[(C1)]\label{def2} $A = \bigoplus_{g \in G} R u_g$.
	\item[(C2)]\label{def3} $\forall g \in G$, $R u_g = u_g R$ and this is a free left $R$-module of rank $1$.
	\item[(C3)]\label{def4} The direct sum $A = \bigoplus_{g \in G} R u_g$ turns $A$ into a $G$-graded ring with $R = A_e$.
\end{description}
We call a ring $A$ fulfilling these properties a \textbf{Pre-Crystalline Graded Ring}.}\\

{\prop \label{def5} With conventions and notation as in Definition \ref{def1}:
\begin{enumerate}
	\item For every $g \in G$, there is a set map $\sigma_g : R \rightarrow R$   defined by: $u_g r = \sigma_g(r)u_g$ for $r \in R$.  The map $\sigma_g$ is in fact a surjective ring morphism.  Moreover, $\sigma_e = \textup{Id}_{R}$.
	\item There is a set map $\alpha : G \times G \rightarrow R$ defined by $u_g u_h = \alpha(g,h)u_{gh}$ for $g,h \in G$.  For any triple $g,h,t \in G$ the following equalities hold:
		\begin{eqnarray}
		\alpha(g,h)\alpha(gh,t)&=&\sigma_g(\alpha(h,t))\alpha(g,ht) \label{def6},\\
		\sigma_g(\sigma_h(r))\alpha(g,h)&=& \alpha(g,h)\sigma_{gh}(r) \label{def7}.
		\end{eqnarray}
	\item $\forall g \in G$ we have the equalities $\alpha(g,e) = \alpha(e,g) = 1$ and $\alpha(g,g^{-1}) = \sigma_g(\alpha(g^{-1},g)).$
\end{enumerate}
}
\begin{flushleft}\textbf{Proof}\end{flushleft} See \cite{NVO6}. $\blok$\\

{\prop Notation as above, the following are equivalent:
\begin{enumerate}
	\item $R$ is $S(G)$-torsionfree.
	\item $A$ is $S(G)$-torsionfree.
	\item $\alpha(g,g^{-1})r=0$ for some $g \in G$ implies $r = 0$.
	\item $\alpha(g,h)r=0$ for some $g,h \in G$ implies $r = 0$.
	\item $R u_g = u_g R$ is also free as a right $R$-module with basis $u_g$ for every $g \in G$.
	\item for every $g \in G$, $\sigma_g$ is bijective hence a ring automorphism of $R$.
\end{enumerate}
}
\begin{flushleft}\textbf{Proof}\end{flushleft} See \cite{NVO6}. $\blok$\\

{\defi Any $G$-graded ring $A$ with properties \textbf{(C1),(C2),(C3)}, and which is $G(S)$-torsionfree is called a \textbf{crystalline graded ring}.  In case $\alpha(g,h) \in Z(R)$, or equivalently $\sigma_{gh}=\sigma_g \sigma_h$, for all $g,h \in G$, then we say that $A$ is \textbf{centrally crystalline}.}\\

We associate the following function to the $2$-cocycle $\alpha$:
\[f_\alpha: G \times G \rightarrow K: (x,s)\mapsto \frac{\alpha(x,s)}{\alpha(xsx^{-1},x)}.\]
We will drop the subscript $\alpha$ if there is no confusion possible.  For the remainder we also need the following definitions:

{\defi A \textbf{regular element of $G$ with respect to $\alpha$} is an $x \in G$ such that $\alpha(x,g)=\alpha(g,x)$ for every $g \in C(x):=\{g \in G| gx=xg\}$.  The set of $\alpha$-regular elements in $G$ is denoted by $G_\textup{reg}$.}\\

{\defi An \textbf{$\alpha$-ray class} is defined to be the conjugation class of a regular element with respect to $\alpha$.}\\

{\defi \label{cen1} Given an $\alpha$-ray class, we define \textbf{the $\alpha$-ray class sum} as the sum of the base elements $u_g$ where $g$ is in the $\alpha$-ray class.}\\

\section{The Center}

\subsection{Setting}
We will use the following notation and conventions.
\begin{itemize}
	\item $R$ is a commutative domain.
	\item $K = \Q(R)$ is the field of fractions of $R$.
	\item $G$ is a finite group.
	\item $A \CGR R$ is a crystalline graded ring.
	\item $W = \textup{Ker} \sigma$. 
\end{itemize}

\subsection{The Center}

For the case of a twisted group ring, see for example \cite{NVO5}.  From now we consider $A= R \CGR G$ where $R$ is a commutative domain and $G$ finite.  Let $W = \textup{Ker}\sigma$.  If we write an inverse, we mean the inverse in the field of fractions $K$ of $R$.

\newpage

{\lem \label{cen2} We have the following formulas $\forall x,g \in G$:
\begin{enumerate}
	\item In $K \CGR G : u_x^{-1} = u_{x^{-1}}\alpha^{-1}(x,x^{-1}) = \alpha^{-1}(x^{-1}, x)u_{x^{-1}}$.
	\item $\alpha(gx, x)\sigma_{xgx^{-1}}(\alpha^{-1}(x, x^{-1}))=\alpha^{-1}(xgx^{-1}, x)$.
\end{enumerate}}
\begin{flushleft}\textbf{Proof}\end{flushleft}
\begin{enumerate}
	\item $u_x u_{x^{-1}}\alpha^{-1}(x,x^{-1}) = \alpha(x,x^{-1})u_e \alpha^{-1}(x,x^{-1}) =1$ and $u_{x^{-1}}\alpha^{-1}(x,x^{-1}) u_x = \sigma_{x^{-1}}(\alpha^{-1}(x,x^{-1}))\alpha(x^{-1},x)=1$.
	\item Use the $2$-cocycle relation (\ref{def6}) for $(xg, x^{-1},x)$. $\blok$
\end{enumerate} 

{\prop \label{cen3} With notation as above we have
\[
\sum\limits_{s \in G} {r_s u_s  \in Z(A)}  \Leftrightarrow \left\{ {\begin{array}{*{20}c}
   {r_s  = 0} \hfill & {\forall s \notin W} \hfill  \\
   {\sigma _x (r_s )\alpha (x,s) = r_{xsx^{ - 1} }\alpha (xsx^{ - 1} ,x) } \hfill & {\forall x \in G, s\in W}. \hfill  \\

 \end{array} } \right.
\]\ \\}
\textbf{Proof}
\begin{itemize}
	\item Let $\sum r_s u_s \in Z(A)$, then $\forall r \in R$:
	\begin{eqnarray*}
	\left(\sum_{s \in G}r_s u_s\right)r = r\left(\sum_{s \in G}r_s u_s\right)&\Leftrightarrow& \sum_{s \in G}r_s \sigma_s(r)u_s = \sum_{s \in G}r_s r u_s\\
	\ &\Leftrightarrow& \sigma_s(r)=r \ \ \ \forall s\in G \textrm{ with } r_s \neq 0\\
	\ &\Leftrightarrow& r_s = 0\ \ \ \forall s\notin W.
	\end{eqnarray*}
	\item Let $x \in G$ and $s \in W$.  We will use the formulas in Lemma \ref{cen2}:
	\begin{eqnarray*}
	u_x r_s u_s u_x^{-1} &=& \sigma_x(r_s)u_x u_s u_{x^{-1}}\alpha^{-1}(x, x^{-1})\\
	\ &=& \sigma_x(r_s)\alpha(x,s)\alpha(xs, x^{-1})u_{xsx^{-1}}\alpha^{-1}(x, x^{-1})\\
	\ &=& \sigma_x(r_s)\alpha(x,s)\alpha(xs, x^{-1})\sigma_{xsx^{-1}}(\alpha^{-1}(x, x^{-1}))u_{xsx^{-1}}\\
	\ &=& \sigma_x(r_s)\alpha(x,s)\alpha(xs, x^{-1})\sigma_{xs}(\alpha( x^{-1},x))^{-1}u_{xsx^{-1}}\\
	\ &=& \sigma_x(r_s)\alpha(x,s)\alpha(xsx^{- 1}, x)u_{xsx^{-1}}.
	\end{eqnarray*}
\end{itemize}
And this finishes the proof. $\blok$\\

If $g \in G$ is a degree of a nonzero element in $Z(A)$, say $0 \neq r_g u_g \in Z(A)$, then $\sigma_x(r_g)\alpha(x,g) = r_g \alpha(g,x)$ for $x \in G$ follows, or $\sigma_x(r_g) = r_g f(g,x)$.  We now calculate in a straightforward way:
\begin{eqnarray*}
r_g f(g,xy) &=& \sigma_{xy}(r_g)\\
\ &=&\sigma_x(\sigma_y(r_g))\\
\ &=& \sigma_x(r_g f(g,y))\\
\ &=& \sigma_x(r_g)\sigma_x(f(g,y))\\
\ &=& r_g f(g,x) \sigma_x(f(g,y)).
\end{eqnarray*}
Hence we arrive at (since $R$ is a domain):
\[f(g,xy)=f(g,x)\sigma_x(f(g,y)).\]
A $g\in G$ corresponding to a nonzero $r_g u_g \in Z(A)$ corresponds therefore to a crossed homomorphism $f(g, -): G \rightarrow K$, i.e. to an element of $H^1(G, K)$ where $K$ is a $G$-module via the extension of the $G$-action on $R$ to the fraction field.  The restriction of $f(g,-)$ to $W \subset G$ defines a multiplicative map $W \rightarrow K$.  This situation appears often when $G$ is Abelian, the case we first treat in some detail hereafter.

\subsection{$G$ Abelian}

From now we consider $A= D \CGR G$ where $D$ is a Dedekind domain and $G$ Abelian finite.  Let $K$ be the field of fractions of $D$.  In this section we will prove that the center of $A$ is crystalline graded itself in certain cases.  It is easy to see that the center is graded since $G$ is Abelian.  With $G$ Abelian, the second statement in Proposition \ref{cen3} becomes
\begin{equation}\sigma_x(d_s)\alpha(x,s)= d_s \alpha(s,x) \ ,\ \ \forall x \in G, s \in W.\label{cen4}\end{equation}
And the function $f := f_\alpha$ becomes
\[f: G \times G \rightarrow K : (g,h)\mapsto \frac{\alpha(g,h)}{\alpha(h,g)}.\]
So rewriting (\ref{cen4}):
\begin{equation}\sigma_x(d_s)=d_s f(s,x) \ ,\ \ \forall x \in G, s \in W.\label{cen5}\end{equation}
Consider $B = D \CGR W$.  The first statement in (\ref{cen3}) gives us $Z(A)\subset Z(B)$.  Since the action on $Z(B)$ is trivial, we find that, using (\ref{cen5})
\[Z(B)=\sum_{s \in W_{\textup{reg}}} D u_s,\]
where $W_{\textup{reg}} = \{s \in W | \alpha(s,x)=\alpha(x,s), \forall x \in W\}$.  Fix a $d_s u_s \in Z(A)$, so $s \in K_{\textup{reg}}$ (looking at these elements is sufficient since the center is graded).  We have two cases: either $s \in G_{\textup{reg}}$, or $s \notin G_{\textup{reg}}$.\\
If $s \in G_{\textup{reg}}$, we have that
\[\sigma_x(d_s) = d_s\ ,\ \ \forall x \in G \Rightarrow d_s \in D^G,\]
where we define $D^G = \{d \in D | \sigma_x(d)=d \ ,\forall x \in G\}$.\\
Now consider the case that $s \notin G_{\textup{reg}}$ and define $\forall x$ with $\alpha(s,x)\neq\alpha(x,s)$ the set $I_{s,x}$:
\[I_{s,x}=\{d\in D | \sigma_x(d) = d f(s,x)\},\]
and define $I_s$ as the intersection of all $I_{s,x}$, $x \in G$ with $\alpha(s,x)\neq \alpha(x,s)$.  So now we have
\[Z(A) = \sum_{s \in W \cap G_{\textup{reg}}}D^G u_s + \sum_{s \in W_{\textup{reg}}\backslash G_{\textup{reg}}}I_s u_s.\]

{\prop \label{cen6}$I_s$ as defined above is a finitely generated $D^G$-bimodule of rank $1$.  It is not multiplicatively closed.}\\
\textbf{Proof}The fact that $I_s$ is a $D^G$-bimodule is easily checked.  Not multiplicatively closed follows from $\sigma_x(mn)= \sigma_x(m)\sigma_x(n)=mnf(x,s)^2 \neq mnf(x,s)$ since $0$ and $1$ are the only idempotents.\\
Since $D$ is finitely generated over $D^G$ and $D^G$ is a Dedekind domain ($G$ finite), we see that $I_s$ has finite rank.  The rank is determined as follows.  Take an $x$ such that $\alpha(s,x)\neq \alpha(x,s)$ and take $u,v \in I_s$:
\[\frac{\sigma_x(u)}{\sigma_x(v)}=\frac{uf(s,x)}{vf(s,x)}= \frac{u}{v},\]
and this entails that $u/v \in K^G$ so the rank of $I_s$ is $1$. $\blok$\\

{\st $Z(A)$ is crystalline graded over $D^G$ when $D$ is a principal ideal domain.}\\
\textbf{Proof} This actually follows almost immediately from Proposition \ref{cen6}.  Since $I_s$ is a torsionfree module over a principal ideal domain, it is free.  Since it has rank $1$ we find that $I_s = D^G d_s$ for some $d_s \in D$ and so our center becomes when we set $v_s = d_s u_s$
\[Z(A) = \sum_{s \in W \cap G_{\textup{reg}}}D^G u_s + \sum_{s \in W_{\textup{reg}}\backslash G_{\textup{reg}}}D^G v_s.\] $\blok$\\

We can now address to the question when the center is trivial, e.g. $Z(A)=D$ or $Z(A)=D^G$.  This means that the only component that can appear is the component corresponding to $u_e=1$, and so our center will be $D^G$.  There are a number of possibilities when this may happen:
\begin{enumerate}
	\item $W_\textup{reg} = \{e\}$.
	\item $W_\textup{reg} \cap G_\textup{reg} = \{e\}$ and $I_s = \{0\} \ \ \forall s \in W_\textup{reg} \backslash G_\textup{reg}$.
\end{enumerate}
We are not quite interested in the first condition, namely that $W$ has no regular elements, but we are interested in the second condition:
\begin{itemize}
	\item $\{e\}=W_\textup{reg}\cap G_\textup{reg} = G_\textup{reg} \cap W$.  So no $G$-regular elements distinct from $e$ may be in $W$.
	\item $I_s = \{0\} \ \ \forall s \in W_\textup{reg} \backslash G_\textup{reg}$.  This means that there is no solution to
	\[\frac{\sigma_x(d)}{d} = f(x,s) \ \ \forall x \in G \textup{ with } \alpha(x,s)\neq\alpha(s,x).\]
\end{itemize}

\subsection{$G$ Not Abelian}

When $G$ is not Abelian, we get a few extra problems.  For one, the center will not be graded.  Furthermore, we will need a lot more criteria to form the center.  For this section, we will set $R=D$ to be a Dedekind domain.  We set
\[A = D \CGR G,\]
\[B = D \CGR W,\]
where $\sigma$ and $\alpha$ are as usual, $K = \textup{Ker}\sigma$.  
We have the following theorem from \cite{CVO}:
{\st Let $R$ be a domain and let $G$ be a finite group.  If a $2$-cocycle $\alpha \in Z^2(G, R)$ has the property that the corresponding function $f_\alpha(x,s)=1$, $\forall \alpha$-regular $s \in G$ and all $x \in G$, then the $\alpha$-ray class sums form an $R$-basis for the center $Z(R *_\alpha G)$ of the generalized crossed product $R *_\alpha G := R \mathop \diamondsuit \limits_{\textup{Id} ,\alpha }G$.}\\

We can use this theorem on $B$ since the twist generated by $\sigma$ is trivial for $W$.  We will call the property that
\begin{equation}f_\alpha(x,s)=1, \ \ \forall s \in G_{\textup{reg}},\ \ \forall x \in G,\label{cen7}\end{equation}
the universal regularity condition \textbf{URC} of $\alpha$ for $G$.  (Remark: for now, $G_{\textup{reg}}$ has only meaning when the twist $\sigma$ is trivial.)  Now denote $C_1, \ldots, C_t$ the ray classes of $W$ and $v_i:=v_{C_i}$ the ray class sum corresponding to $C_i$, $\forall i$.  We find, if $\alpha$ satisfies the \textbf{URC} (\ref{cen7}) for $W$:
\[Z(B)= \sum_{i=1}^{t}D v_i.\]
{\lem With notations as above, we have the following formula $\forall x \in G$, $\forall i = 1,\ldots, t$:
\[u_x v_i u_x^{-1} = \sum_{g \in C_i} f(x,g)u_{xgx^{-1}}.\]}
\begin{flushleft}\textbf{Proof}\end{flushleft}
We use the formulas in Lemma \ref{cen2}
	\begin{eqnarray*}
     u_x v_i u_x^{-1} &=& u_x v_i u_{x^{-1}}\alpha^{-1}(x,x^{-1})\\
     \ &=& u_x \left(\sum_{g \in C_i}u_g\right)u_{x^{-1}}\alpha^{-1}(x,x^{-1})\\
     \ &=& \left(\sum_{g \in C_i}\alpha(x,g)u_{xg}\right)u_{x^{-1}}\alpha^{-1}(x,x^{-1})\\
     \ &=& \sum_{g \in C_i}\alpha(x,g)\alpha(xg,x^{-1})\sigma_{xgx^{-1}}(\alpha^{-1}(x,x^{-1}))u_{xgx^{-1}}\\
     \ &=& \sum_{g \in C_i}\alpha(x,g)\alpha^{-1}(xgx^{-1}, x)u_{xgx^{-1}}\\
     \ &=& \sum_{g \in C_i} f(x,g)u_{xgx^{-1}}.
  \end{eqnarray*}  $\blok$\\

For the remainder of this section, let us assume that $\alpha$ satisfies the \textbf{URC} (\ref{cen7}) for $W$.  We can now define
\[W_{\textup{reg}}=\{g \in W | f(x,g)=1 \ \ \forall x \in W\},\]
\[G_{\textup{reg}}=\{g \in G | f(x,g)=1 \ \ \forall x \in G\}.\]
Like in the previous section, we can use Proposition \ref{cen3} so we restrict to $Z(B)$.  \\

{\lem \label{cen8} If the $2$-cocycle $\alpha$ satisfies the \textbf{URC} (\ref{cen7}), we have that $\forall x \in G$, $\forall s \in W_\textup{reg}$, $xsx^{-1} \in W_\textup{reg}$.}\\
\textbf{Proof} Since $\alpha$ satisfies the \textbf{URC} (\ref{cen7}) we only have to check that for $g \in C_W(xsx^{-1})$ we have $u_g u_{xsx^{-1}} = u_{xsx^{-1}} u_g$.  (One easily checks that this is equivalent with regularity.)  We see that $x^{-1}gx \in C_W(s)$ and thus since $s \in W_\textup{reg}$ (we use the formulas in Lemma \ref{cen2})
\begin{eqnarray*}
\ &\ & u_s u_{x^{-1}g x} = u_{x^{-1}g x}u_s\\
&\Rightarrow& u_s \alpha^{-1}(x^{-1}g,x)\alpha^{-1}(x^{-1},g)u_{x^{-1}}u_g u_x \\
&\ &\ \ \ \ \ \ \ \ \ \ \ \ \ \ \ \ \ = \alpha^{-1}(x^{-1}g,x)\alpha^{-1}(x^{-1},g)u_{x^{-1}}u_g u_x u_s \\
&\Rightarrow& u_s u_{x^{-1}}u_g u_x = u_{x^{-1}}u_g u_x u_s\\
&\Rightarrow& u_{x^{-1}}^{-1}u_s u_{x^{-1}}u_g= u_g u_x u_s u_x^{-1}\\
&\Rightarrow& \alpha^{-1}(x,x^{-1})u_x u_s u_{x^{-1}}u_g = u_g u_x u_s u_{x^{-1}}\alpha^{-1}(x,x^{-1})\\
&\Rightarrow& u_x u_s u_{x^{-1}}u_g = u_g u_x u_s u_{x^{-1}}\\
&\Rightarrow& u_{xsx^{-1}}u_g = u_g u_{xsx^{-1}}.
\end{eqnarray*}
And so $xsx^{-1} \in W_\textup{reg}$. $\blok$\\

Take an element $y =\sum_{i = 1}^t d_i v_i \in Z(B)$, $d_i \in D$.  We have to take a full sum now, since maybe the direct sum we have in $Z(B)$ does not hold in $Z(A)$.  We calculate
\begin{equation}u_x \left(\sum_{i = 1}^t d_i v_i\right) u_x^{-1}= \sum_{i=1}^t \sigma_x(d_i)\sum_{g\in C_i}f(x,g)u_{xgx^{-1}}.\label{cen9}\end{equation}
So for $\sum_{i=1}^t d_i v_i$ to be an element of $Z(A)$ we must have that this equation (\ref{cen9}) equals $\sum_{i=1}^t d_i v_i$, $\forall x \in G$.  In other words, this means that for all appearing ray classes $C_i$, we have that $\forall g \in C_i$, $\forall x \in G$ that $xgx^{-1} \in W_\textup{reg}$.  This is true according to Lemma \ref{cen8}.  We get the following condition on the coefficients $d_i$ of $y$ to be in the center of $A$, $\forall x \in G, \forall i = 1,\ldots, t,\forall g \in C_i$:
\begin{equation}\sigma_x(d_i)f(x,g) = d_j\ \ ,\ \ \textrm{ where } j \textrm{ is determined by } xgx^{-1} \in C_j. \label{cen10}\end{equation}
Let $d_i \neq 0$ for some $i$.  We find, looking at equation (\ref{cen10}), that all coefficients $d_j$ corresponding to the ray classes $C_j$ will be nonzero, where $j$ is determined by $xgx^{-1} \in C_j$ for some $x\in G$, some $g \in C_i$.  Define for ease of notation the saturation class $\Gamma_i$, $\forall i = 1, \ldots, t$ as the set consisting of all ray classes $C$ (for $W$) such that if you take $x \in G$, $C$ in $\Gamma_i$ and $g \in C$, there is a $C' \in \Gamma_i$ with $xgx^{-1} \in C'$.  So this means the general expression for an element $y \in Z(A)$ may be written as:
\[y = \sum_{C \in \Gamma_i}d_C v_C\ \ \ \ \textrm{ for some }i.\]
Let $\bar x = \bar z \in G/W$ or equivalently $\exists k \in W: z =kx$.\\

{\lem Let $g,k \in W$, and $x \in G$.  Then
\[f(kx,g)=f(x,g)f(k, xgx^{-1}),\]
i.e. if $g \in W_\textup{reg}$ then $f(kx,g)=f(x,g)$.}\\
\textbf{Proof}
\begin{eqnarray*}
f(kx,g) &=& \frac{\alpha(kx,g)}{\alpha(kxgx^{-1}k^{-1}, kx)} = \frac{\alpha(kx,g)u_{kxg}}{\alpha(kxgx^{-1}k^{-1}, kx)u_{kxg}}=\frac{u_{kx}u_{g}}{u_{kxgx^{-1}k^{-1}}u_{kx}}\\
\ &=& \frac{\alpha^{-1}(k,x)u_k u_x u_g}{\sigma_{kxgx^{-1}k^{-1}}(\alpha^{-1}(k,x))u_{kxgx^{-1}k^{-1}}u_{k}u_{x}}\\
\ &=& \frac{\alpha(x,g)u_ku_{xg}}{\alpha(kxgx^{-1}k^{-1},k)\alpha^{-1}(k,xgx^{-1})u_k u_{xgx^{-1}}u_x}\\
\ &=& \frac{\alpha(x,g)u_k u_{xg}}{\alpha(xgx^{-1}, x)\alpha(kxgx^{-1}k^{-1},k)\alpha^{-1}(k, xgx^{-1})u_k u_{xg}}\\
\ &=& f(x,g)f(k,xgx^{-1}).
\end{eqnarray*}
$\blok$\\

{\gev If $\forall g \in C$, $C$ ray class in $W$, $\forall x,z \in G$ such that $\bar x = \bar z$ in $G/W$, then
\[f(x,g) = f(z,g).\]}\\

In particular, for $g \in C$, ray class in $W$ we have that
\[F_g : G/W \rightarrow K : \bar x \mapsto f(x,g)\]
is well defined.\\

Since the condition (\ref{cen10}) is not as beautiful as the equation in the Abelian case (\ref{cen5}), we will not describe the center fully.  Instead we want to find conditions on the $2$-cocycle $\alpha$ and twist $\sigma$ to have a trivial center, e.g. $Z(A)=D$ or $Z(A)=D^G$.  This for example happens if $W = \{e\}$.  We can also describe the center when $W = G$, since there will be no twist and we can use the theory of twisted group rings.  We can first assume that the center is of the following form:
\[Z(A) = \bigoplus_{C \textup{ a certain ray class in }W}d_C v_C.\]\
This actually implies that a general element $y \in Z(A)$ is of the form $d_C v_C$ for a certain ray class $C$ in $W$.  Not all ray classes $C$ will be present though.  Looking at equation (\ref{cen9}), for $y$ to be in $Z(A)$, we need to have if $d_C \neq 0$ then $\forall x \in G: xCx^{-1} = C$.  In group theoretic words this means that the normalizer $N_G(C) = \{x \in G | xCx^{-1} = C\}=G$.  So if $N_G(C)\neq G$, the component corresponding to $C$ will be $0$.  So for now fix $C$ with $N_G(C)=G$ and set $u_x (d_C v_
C) u_x^{-1} = d_C v_C$ then we find
\[\sigma_x(d_C)f(x,g) = d_C \ \ \ \forall g \in C.\]
If $d_C \neq 0$, then $f(x,g) = f(x,h)$ $\forall g,h \in C$, $\forall x \in G$.  This in particular means that if $d_C \neq 0$, $f(x,g)$ with $g \in C$ is fully determined by one representative, say $s \in C$.  We will call such $C$ constant under $f$.  We define
\[\Delta = \{C | C\textrm{ is a ray class in }W \textrm{ with } C \textrm{ constant under }f\}.\]
So from now on, fix $C \in \Delta$.  We can now define $f(x, C)$ as $f(x, s)$ for some $s \in C$.  Like in the Abelian case, we have two possibilities: $f(x, C)=1$, $\forall x \in G$ (we say $C \in G_{\textup{reg}}$) or $\exists x \in G: f(x,C)\neq 1$ (we say $C \notin G_{\textup{reg}}$).  If $C \in G_{\textup{reg}}$ we see that $\sigma_x(d_C)=d_C$, $\forall x \in G$ or that $d_C \in D^G$.  If $C \notin G_{\textup{reg}}$ we define for each $x \in G$ with $f(x,C)\neq 0$:
\[I_{C,x}=\{d\in D | \sigma_x(d) = d \cdot f^{-1}(x,C)\},\]
and define $I_{C}$ to be the intersection of all $I_{C,x}$, $x\in G, f(x,C)\neq 0$.  Similar to the Abelian case, our center now becomes:
\[Z(D \CGR G) = \sum_{C \in G_{\textup{reg}} \cap \Delta}D^G v_C + \sum_{C \in \Delta \backslash G_{\textup{reg}}}I_{C}v_C.\]
One can prove in exactly the same way as in the Abelian case:\\

{\prop $I_{C}$ as defined above is a finitely generated $D^G$-bimodule of rank $1$.  It is not multiplicatively closed.}\\
In order to have a trivial center, there are three possibilities:
\begin{enumerate}
  \item $N_G(C)=G \Leftrightarrow C = \{e\}$.
	\item $\Delta = \{C_e\}$.
	\item $G_\textup{reg} = \{C_e\}$ and $I_C = \{0\}$, $\forall C \in \Delta \backslash G_{\textup{reg}}$.
\end{enumerate}
 
{\opm Notice how similar these conditions are to having a trivial center in the Abelian case.}\\

Now let us return to the bimodule considerations following Definition \ref{cen1} now restricting to the situation where $R =D$ is a Dedekind domain.  Since in a Dedekind domain any nonzero ideal is invertible, it follows that $I(g,g^{-1})$ as well as $I(g^{-1},g)$ are invertible and then $A_g A_{g^{-1}}=I(g,g^{-1})$ entails $A_g(A_{g^{-1}}I(g,g^{-1})^{-1}) = D$.  One easily calculates:
\begin{eqnarray*}
(A_{g^{-1}}I(g, g^{-1})^{-1})A_g &=& \sigma_{g^{-1}}(I(g,g^{-1})^{-1})A_{g^{-1}}A_g\\
\ &=& \sigma_{g^{-1}}(I(g, g^{-1})^{-1})I(g^{-1}, g).
\end{eqnarray*}
From $I(g, g^{-1})=\sigma_g(I(g^{-1},g))$ it then follows that the foregoing reduces to
\[A_{g^{-1}}I(g,g^{-1})^{-1}A_g = I(g^{-1},g)^{-1}I(g^{-1},g) = D.\]
Consequently each $A_g$ is an invertible $D$-bimodule with inverse $A_{g^{-1}}I(g,g^{-1})^{-1}$.  Clearly $A$ need not be strongly graded but the equation
\begin{equation}I(X,Y)I(XY,Z)= \sigma_X(I(Y,Z))I(X,YZ),\label{bla}\end{equation}
now gives rise to a map
\[\phi:G \rightarrow \textup{Pic}(D): g \mapsto [A_g],\]
where $[A_g]$ denotes the class of the invertible $D$-bimodule $A_g$ in the class group (Picard group) $\textup{Pic}(D)$, together with a factor system 
\[I: G \times G \rightarrow \textup{Pic}(D):(g,h) \mapsto [I(g,h)].\]
The above equation (\ref{bla}) $I$ defines a $2$-cocycle.  Recall, $\textup{Pic}(D)$ is an Abelian group with respect to the operation induced by the tensor product, it may also be seen as the free Abelian group $\textup{Div}(D)$ generated by the prime ideals of $D$ modulo the subgroup generated by the principal ideals, in particular for a P.I.D. the Picard group is trivial.
Hence $A$ may be viewed as a generalized crossed product
\[A = \bigoplus_{g\in G} A_g,\]
defined by $\phi : G \rightarrow \textup{Pic}(D)$ and a $2$-cocycle $I$ in $\textup{Pic}(D)$, $I : G\times G \rightarrow \textup{Pic}(D)$.  Let us call this type of crossed product a \textbf{class crossed product}.  The notion of crystalline graded ring may be extended to include the cases where the $A_g$ are not necessarily free of rank $1$ but represent elements of the class group $\textup{Pic}(D)$ as above.  We refer to such graded rings as being \textbf{arithmetically crystalline graded}.  The following theorem is now just a rephrasing of the results we obtained:

{\st \label{cen11}Let $A$ be crystalline graded by an Abelian group $G$ over a Dedekind domain $D$, then $Z(A)$ is arithmetically crystalline graded over a subgroup of $G$.}\\

The statement generalizes to $A$ which are themselves arithmetically crystalline graded.  A general definition, preceding those of generalized Weyl algebras and generalized crossed products was used by second author in some observations concerning strongly graded rings.  The so-called $\delta$-strongly graded rings were characterized by $A_g A_{g^{-1}} = \delta_g$ being an invertible $R$-bimodule for each $g\in G$.  Of course, if $R$ is not a Dedekind domain and in particular if it is not even commutative then $\textup{Pic}(R)$ is somewhat more complex to deal with, for example a $[M] \in \textup{Pic}(R)$ canonically defines an automorphism $\sigma_M$ of $Z(R)$ but not necessarily of $R$ (note that in the situation of Theorem \ref{cen11} the $\sigma_{A_g}$ are exactly the $\sigma_g$ defined on $R$!).  More detail about $\textup{Pic}(R)$ in the greater generality can be found in the book by H. Bass (K-theory), cf. \cite{B}.  For us this will be useful in further work on the algebraic structure of general crystalline or arithmetically crystalline graded rings e.g. over Artinian algebras etc.

\end{document}